\documentclass[final,leqno]{siamltex}

\usepackage{stmaryrd}
\usepackage{amsmath}
\usepackage{amssymb}
\usepackage{multirow}
\usepackage{cite}
\usepackage{tabularx}
\usepackage{color}

\usepackage{comment}
\usepackage{tabularx}

\usepackage{algorithm}
\usepackage{algorithmic}

\newtheorem{remark}[theorem]{Remark}

\newcommand\bx{\mathbf{x}}
\newcommand\bn{\mathbf{n}}
\newcommand\bm{\mathbf{m}}

\newcommand\bk{\mathbf{k}}

\newcommand\by{\mathbf{y}}
\newcommand\bR{\mathbb{R}}
\newcommand\Udip{U_{\text{dip}}}

\newcommand\hf{\hat{f}}
\newcommand\wu{\widehat{U}}
\newcommand\wrho{\widehat{\rho}}

\title{Fast and Accurate Evaluation of  Nonlocal
Coulomb and  Dipole-Dipole  Interactions via the Nonuniform FFT}
\author{Shidong Jiang\thanks{Department of Mathematical Sciences,
New Jersey Institute of Technology, Newark, New Jersey 07102
({\tt shidong.jiang@njit.edu}). This
research was supported by the National Science Foundation under grant CCF-0905395. }
\and Leslie Greengard\thanks{Courant Institute of Mathematical
Sciences, New York University, New York, New York 10012
({\tt greengard@courant.nyu.edu}). This research was supported in part by the
U.S. Department of Energy under contract DEFG0288ER25053.}
\and Weizhu Bao\thanks{Department of Mathematics and Center for
Computational Science and Engineering, National
University of Singapore, Singapore 119076, Singapore
({\tt matbaowz@nus.edu.sg}). This
research was supported by the Singapore A*STAR SERC  PSF-Grant 1321202067. }}

\begin{document}

\maketitle

\begin{abstract}
We present a fast and accurate algorithm for the evaluation of
nonlocal (long-range) Coulomb and dipole-dipole interactions in free space.
The governing potential is simply
the convolution of an interaction kernel $U(\bx)$ and
a density function $\rho(\bx)=|\psi(\bx)|^2$, for some
complex-valued wave function $\psi(\bx)$, permitting the
formal use of Fourier methods. These are hampered by the fact that
the Fourier transform of the interaction kernel
$\widehat{U}(\bk)$ has a singularity at the origin $\bk={\bf 0}$
in Fourier (phase) space. Thus, accuracy is lost when using a uniform
Cartesian grid in $\bk$ which would otherwise
permit the use of the FFT for evaluating the convolution.
Here, we make use of a high-order discretization of the Fourier integral,
accelerated by the nonuniform fast Fourier transform (NUFFT).
By adopting spherical and polar phase-space discretizations
in three and two dimensions, respectively, the
singularity in $\hat{U}(\bk)$ at the origin is canceled, so that
only a modest number of degrees of freedom are required to evaluate
the Fourier integral, assuming
that the density function $\rho(\bx)$ is smooth and decays sufficiently fast
as $\bx \rightarrow \infty$. More precisely, the calculation requires
$O(N\log N)$ operations, where $N$ is the total number of
discretization points in the computational domain.
Numerical examples are presented to demonstrate
the performance of the algorithm.
\end{abstract}

\begin{keywords}
Coulomb interaction,  dipole-dipole interaction, interaction energy,
nonuniform FFT, nonlocal, Poisson equation.
\end{keywords}

\begin{AMS}
33C10, 33F05, 44A35, 65R10, 65T50, 81Q40
\end{AMS}

\pagestyle{myheadings}
\thispagestyle{plain}
\markboth{S. Jiang, L. Greengard, and W. Bao}
{Fast Evaluation for Nonlocal Interactions via NUFFT}

\section{Introduction} \label{introduction}

Nonlocal (long-range) interactions are encountered in modeling a variety of problems
from quantum physics and chemistry to materials science and biology.
A typical example is the Coulomb interaction in the nonlinear Schr\"{o}dinger equation
(or Schr\"{o}dinger-Poisson system in three dimensions (3D))
as a ``mean field limit'' for $N$-electrons, assuming binary
Coulomb interactions \cite{BEGMY,BGM,HF} and the Kohn-Sham equation of
density functional theory (DFT) for electronic structure calculations
in materials simulation and design \cite{Kohn,HF,Parr,Sheehy,Sheehy1}.
Dipole-dipole interactions arise in quantum chemistry \cite{Levitt,Hakn}, in dipolar
Bose-Einstein condensation (BEC) \cite{dip1,dip2,gpe3,gpe4,gpe5,gpe6,BC7,bao1,bao2,Lahaye},
 in dipolar Fermi gases \cite{Miyakawa}, and in dipole-dipole interacting Rydberg molecules
\cite{Kiffner1,Kiffner2,Kiffner3}.

In physical space, the interaction kernel is both long-range and singular at the
origin, requiring both accurate quadrature techniques and suitable fast algorithms.
When the density function is smooth, however, it is often more convenient to use Fourier
methods since the frequency content is well-controlled. Unfortunately, the Fourier
transform of the interaction kernel is singular at the origin of Fourier (phase) space
as well, resulting in significant numerical burdens and challenges
\cite{bao1,BC7,BMS,nufft1,Blakie,fracpoisson1,Tick,SPMCompare}.

In this paper, we present a fast and accurate algorithm for the numerical evaluation
of the interaction potential \cite{BEGMY,HF,Kohn,Parr,Sheehy,gpe3,gpe4,gpe5,gpe6}:
\begin{equation}\label{gpe}
u(\bx)=(U\ast\rho)(\bx):=\int_{\mathbb{R}^d}U(\bx-\by)\rho(\by)\,d\by,
\qquad \bx \in \mathbb{R}^d, \qquad d=3,2,
\end{equation}
and its related interaction energy \cite{BEGMY,Kohn,HF,Parr,Sheehy,gpe3,gpe4,gpe5,gpe6}
\begin{equation}\label{gpe2}
E(\rho):=\frac{\lambda}{2}\int_{\mathbb{R}^d}u(\bx)\rho(\bx)\,d\bx=
\frac{\lambda}{2}\int_{\mathbb{R}^d\times\mathbb{R}^d}\rho(\bx)U(\bx-\by)\rho(\by)\,d\by d\bx,
\end{equation}
where $U(\bx)$ is a nonlocal (long-range) interaction kernel and
$\rho(\bx)=|\psi(\bx)|^2$ is a density function derived from a complex-valued
wave function $\psi(\bx)$. Here, $\lambda$ is a dimensionless interaction constant,
and $\ast$ denotes the convolution operator. In most applications,
the density function $\rho$ is smooth and very rapidly decaying
\cite{Kohn,HF,Parr,Sheehy,gpe3,gpe4,gpe5,bao2}, so that it can be viewed as having
compact support to a prescribed precision $\varepsilon$.
We focus our attention on the following Coulomb and dipole-dipole interactions:

\vspace{.2in}

\begin{enumerate}
\item
{\em Coulomb interactions in 3D \cite{BEGMY,BGM,HF,Kohn,HF,Parr,Sheehy,Bao8}.}
The interaction
kernel and its Fourier transform are given by the formulas
\begin{equation}\label{Cou3D}
U_{\rm Cou}(\bx)=\frac{1}{4\pi|\bx|} \quad \Longleftrightarrow \quad \widehat{U}_{\rm Cou}(\bk)=\frac{1}{\|\bk\|^2},
\qquad \bx,\bk\in \mathbb{R}^3.
\end{equation}
In certain settings, the density is limited to two dimensions and one seeks
the Coulomb potential in that plane alone. This arises in various problems
of surface physics \cite{ReductionNanowire,fracpoisson1,cichoki,Bao8}, and the
governing potential is obtained in two dimensions by dimension reduction from
three dimensions under an anisotropic potential. This is well-known to yield:
\begin{equation}\label{Cou2D}
U^{(2.5)}_{\rm Cou}(\bx)=\frac{1}{2\pi|\bx|} \quad \Longleftrightarrow \quad
\widehat{U}^{(2.5)}_{\rm Cou}(\bk)=\frac{1}{\|\bk\|},
\qquad \bx,\bk\in \mathbb{R}^2.
\end{equation}
The superscript $(2.5)$ is intended to denote that the sources line in
a two-dimensional space but that the physical interaction is that of the
ambient three-dimensional space.
Here $\widehat{f}(\bk)$ is the Fourier transform of a function $f(\bx)$
defined by the formula $ \widehat{f}(\bk)=
\int_{\mathbb{R}^d} f(\bx)\; e^{-i\,\bk\cdot \bx}\;d\bx$ for $\bx,\bk\in \mathbb{R}^d$.

\item
{\em Dipole-dipole interactions with the same dipole orientation in
3D \cite{gpe3,gpe4,gpe5,gpe6,bao1,bao2,BC7,bao3,gpe2,Levitt,Odell0,ODell,Tick,Lahaye}.}
The interaction kernel is given by the formula
\begin{equation}\label{dipkernel}
\Udip(\bx)=\frac{3}{4\pi}\frac{1-3(\bx\cdot \bn)^2/|\bx|^2}{|\bx|^3}
= -\delta(\bx)-3\,\partial_{\bn\bn}\left(\frac{1}{4\pi|\bx|}\right),
\qquad \bx \in \bR^3,
\end{equation}
and its Fourier transform is \cite{bao1,bao2,Odell0,ODell,Xiong,Lahaye}:
\begin{equation}\label{dipkernelF}
\widehat{U}_{\rm dip}(\bk)=-1+\frac{3(\bn\cdot \bk)^2}{\|\bk\|^2},
\qquad \bk \in \bR^3,
\end{equation}
where $\bn=(n_1,n_2,n_3)^T$ is a fixed unit vector representing the dipole orientation,
$\delta$ is the Dirac distribution function, $\partial_\bn=\bn\cdot \nabla$ and
$\partial_{\bn\bn}=\partial_\bn(\partial_\bn)$.
As in the Coulmob case, when the source distribution is two-dimensional, one
reduction from three-dimensions under an
anisotropic potential \cite{bao2,BC7,bao3,bao4,Rosen} yields:
\begin{equation}\label{dipkernel2d}
\Udip^{(2.5)}(\bx)= -\alpha\,
\delta(\bx)-\frac{3}{2}\left(\partial_{\bn_\perp\bn_\perp}-
n_3^2\Delta_\perp\right)\left(\frac{1}{2\pi|\bx|}\right),
\qquad \bx \in \bR^2,
\end{equation}
and its Fourier transform is \cite{bao2,bao3,bao4,Rosen}
\begin{equation}\label{dipkernelF2d}
\widehat{U}^{(2.5)}_{\rm dip}(\bk)=-\alpha+\frac{3\left[(\bn_\perp\cdot \bk)^2-n_3^2\|\bk\|^2\right]}{2\|\bk\|},
\qquad \bk \in \bR^2,
\end{equation}
where $\bn_\perp=(n_1,n_2)^T$, $\partial_{\bn_\perp}=\bn_\perp\cdot \nabla_\perp$,
$\partial_{\bn_\perp\bn_\perp}=\partial_{\bn_\perp}(\partial_{\bn_\perp})$,
$\nabla_\perp=(\partial_x,\partial_y)^T$,
 $\Delta_\perp =\partial_{xx}+\partial_{yy}$, and
$\alpha$ is a fixed real constant \cite{bao2,bao4}.

\item
{\em Dipole-dipole interactions with different dipole
orientations in 3D \cite{Levitt,Hakn,Odell0,ODell,Parker}.}
The interaction kernel is
\begin{eqnarray}\label{dipkerneld}
\qquad \quad \Udip(\bx) &=& \frac{3}{4\pi}\frac{\bm\cdot\bn-3(\bx\cdot \bn)
(\bm\cdot\bx)/|\bx|^2}{|\bx|^3} \\
\nonumber
&=& -(\bm\cdot\bn)\delta(\bx)-3\, \partial_{\bn\bm}\left(\frac{1}{4\pi|\bx|}\right),\ \bx \in \bR^3,
\end{eqnarray}
and its Fourier transform is \cite{Levitt,Hakn,Odell0,ODell,Parker}
\begin{equation}\label{dipkerneldF}
\widehat{U}_{\rm dip}(\bk)=-(\bm\cdot\bn)+\frac{3(\bn\cdot \bk)(\bm\cdot\bk)}{\|\bk\|^2},\qquad \bk \in \bR^3,
\end{equation}
where $\bn=(n_1,n_2,n_3)^T$ and $\bm=(m_1,m_2,m_3)^T$ are two fixed unit vectors
representing the two dipole orientations, $\partial_\bm=\bm\cdot \nabla$ and
$\partial_{\bn\bm}=\partial_\bn(\partial_\bm)=\partial_\bm(\partial_\bn)$.
For a two dimensional distribution, dimension reduction
from 3D under an anisotropic potential \cite{bao2,BC7,bao3,bao4,Rosen} yields:
\begin{equation}\label{dipkernel2dd}
\Udip^{(2.5)}(\bx)= -\alpha\,
\delta(\bx)-\frac{3}{2}\left(\partial_{\bn_\perp\bm_\perp}-
n_3m_3\Delta_\perp\right)\left(\frac{1}{2\pi|\bx|}\right),
\qquad \bx \in \bR^2,
\end{equation}
and its Fourier transform is \cite{bao2,bao3,bao4,Rosen}
\begin{equation}\label{dipkernelF2dd}
\widehat{U}^{(2.5)}_{\rm dip}(\bk)=-\alpha+\frac{3\left[(\bn_\perp\cdot \bk)(\bm_\perp\cdot \bk)-n_3m_3\|\bk\|^2\right]}{2\|\bk\|},
\qquad \bk \in \bR^2,
\end{equation}
where $\bm_\perp=(m_1,m_2)^T$, $\partial_{\bm_\perp}=\bm_\perp\cdot \nabla_\perp$,
$\partial_{\bn_\perp\bm_\perp}=\partial_{\bn_\perp}(\partial_{\bm_\perp})=
\partial_{\bm_\perp}(\partial_{\bn_\perp})$
and $\alpha$ is a fixed real constant \cite{bao2,bao4}.
\end{enumerate}

\vspace{.2in}

\begin{remark}
Note that the second category above is a special case of the third.
It is listed separately because it is simpler and has some important physical
applications.
\end{remark}

\vspace{.2in}

Various numerical methods have been proposed in the literature for evaluating
the interaction potential (\ref{gpe})  and
interaction energy (\ref{gpe2}) using a uniform grid on a bounded
computational domain so as to compute the ground states and dynamics
of problems in quantum physics and chemistry.
By making direct use of the standard uniform FFT
\cite{BMS,Blakie,Odell0,ODell,Parker,Tick,Tikhonenkov,Xiong,gpe4,gpe6},
a phenomoenon know as ``numerical locking" occurs, limiting the achievable precision
\cite{Blakie,Tick,Tikhonenkov,Xiong,bao1,SPMCompare}.
This is due, in essence, to the fact that
$\widehat{U}_{\rm Cou}(\bk)$ and
$\widehat{U}_{\rm dip}(\bk)$ are unbounded at the origin.

As a result, there has been some interest in reformulating the problem of
convolution with the 3D Coulomb kernel
(\ref{Cou3D}) in terms of the governing partial differential equation (the Poisson
equation)
\begin{equation}\label{gpe3D}
-\Delta u_{\rm Cou}(\bx)=\rho(\bx), \qquad \bx \in \mathbb{R}^3,
 \qquad \qquad \lim_{|\bx|\to\infty}u_{\rm Cou}(\bx)=0,
\end{equation}
and convolution with the reduced 2.5D Coulomb kernel
(\ref{Cou2D}) in terms of the fractional partial differential equation
\begin{equation}\label{gpe2D}
(-\Delta)^{1/2} u_{\rm Cou}(\bx)=\rho(\bx), \qquad \bx
\in \mathbb{R}^2, \qquad \qquad \lim_{|\bx|\to\infty}u_{\rm Cou}(\bx)=0.
\end{equation}
The dipole-dipole interaction in 3D (\ref{dipkernel}) can be computed
from the relation \cite{bao1,BC7,Bao8,bao3,bao4,SPMCompare}
\begin{equation}\label{gpeddi}
u(\bx) = -\rho(\bx) +3\partial_{\bn\bn} u_{\rm Cou}(\bx), \qquad \bx \in \mathbb{R}^3.
\end{equation}

There is a substantial literature on solving the
the PDEs (\ref{gpe3D}) and (\ref{gpe2D}),
which we do not seek to review here. We refer the interested reader to
\cite{bao1,BC7,Bao8,bao3,bao4,James,Mccorquodale,SPMCompare} and the
references therein.
We would, however, like to
point out that when the density function $\rho$ has complicated local structure,
an adaptive grid is needed for resolution.
In that setting, Fourier methods are highly inefficient and
the fast multipole method (FMM) or some variant
\cite{fmm3,fmm1,fmm2,fracpoisson1,kifmm2,kifmm5}
can be used for evaluating the nonlocal interaction
directly in physical space in $O(N)$ time, where $N$ is the number of grid points.
In many applications, however, such as the computation of the ground state and
dynamics of BEC \cite{bao1,BC7,bao3,Blakie,bao4,Tick,gpe4,gpe5,gpe6},
one needs to evaluate $u(\bx)$ on an equispaced grid in physical space many times
for different $\rho(\bx)$. This occurs, for example, in
time-splitting spectral methods for computing the dynamics of
the nonlinear Schr\"{o}dinger/Gross-Pitaevskii equations \cite{bao1,BC7,BJM,BJMM}.
In such cases, Fourier methods can be very efficient, easy to implement, and high
order accurate, so long as care is taken in discrettization.

We begin by noting that in Fourier space and the discussion above,
the Coulomb or dipole-dipole interaction
potential \eqref{gpe} is given by
\begin{equation}\label{ifft}
u(\bx)=\frac{1}{(2\pi)^d}\int_{\mathbb{R}^d}  e^{i\,\bk\cdot\bx}\; \widehat{U}(\bk)\; \widehat{\rho}(\bk)\,d\bk,
\qquad \bx\in \mathbb{R}^d, \qquad d=2,3,
\end{equation}
where $\widehat{U}(\bk)$ is given by one of the formulae
\begin{equation} \label{int3dd}
\widehat{U}(\bk)=
\left\{\begin{aligned}
&\frac{1}{\|\bk\|^2}, & \text{3D Coulomb interactions,}\\
&-(\bm\cdot\bn)+\frac{3(\bn\cdot \bk)(\bm\cdot\bk)}{\|\bk\|^2}, & \text{3D dipole-dipole interactions,} \\
&\frac{1}{\|\bk\|}, &  \text{2.5D Coulomb interactions,} \\
&-\alpha+\frac{3\left[(\bn_\perp\cdot \bk)(\bm_\perp\cdot \bk)-n_3m_3\|\bk\|^2\right]}{2\|\bk\|}, &
\text{2.5D dipole-dipole interactions.}
\end{aligned}\right.
\end{equation}

The remainder of this paper is aimed at the construction of a fast and accurate
algorithm for the evaluation
of long-range interactions of the form \eqref{gpe}
as well as the total interaction
energy \eqref{gpe2}. There are three essential ingredients.
First, we truncate the integrals in \eqref{ifft}
at a frequency beyond which the contribution to $\rho$ is negligible.
This is valid because of our assumption that $\rho$ is smooth.
Second, we rewrite
\eqref{ifft} using spherical or polar coordinates in 3D or 2D,
respectively. The Jacobian of this change of variables cancels the singularity
at the origin in Fourier space, permitting the use of simple high order
quadrature rules. More precisely, we achieve superalgebraic convergence by
using the trapezoidal rule in the azimuthal direction and Gauss-Legendre quadrature
in the radial and inclination directions.
Third, we utilize the nonuniform fast Fourier transform (NUFFT) (see, for example,
\cite{nufft1,nufft2,nufft3,nufft4,nufft5,nufft6,nufft7,nufft8}) to accelerate
the calculation of the sums which arise from discretization, which do not correspond
to uniform tensor product grids.
The resulting algorithm is high-order accurate and requires only $O(N\log N)$ work,
where $N$ is the total number of discretization points in physical space.

The paper is organized as follows.
In Section 2, we present a brief review of NUFFT and in Section 3,
we describe the numerical algorithm in detail. The performance of the method
is illustrated with several
numerical examples in Section 4.
Section 5 contains some concluding remarks.

\section{Brief review of the NUFFT}

In this section, we summarize the basic steps of
the NUFFT, to make the discussion reasonably self-contained.

The ordinary FFT computes the discrete Fourier transform (DFT) and its inverse:
\begin{equation}\label{fft}
\begin{aligned}
F(k)&=\sum_{j=0}^{N-1} f(j) e^{-2\pi i k j/N}, & k=0,\cdots,N-1, \\
f(j)&=\frac{1}{N}\sum_{k=0}^{N-1}F(k) e^{2\pi i kj/N },& j=0,\cdots,N-1
\end{aligned}
\end{equation}
in $O(N\log N)$ operations by exploiting the algebraic structure of the
DFT matrix. The points $x_j = 2\pi j/N$ and the frequencies
$k$, however must be equispaced in
both the physical and Fourier domains (see, for example, \cite{fft,strang}).

The purpose of the NUFFT is to remove this restriction,
while maintaining a computational complexity of $O(N\log N)$, where
$N$ denotes the total number of points in both the physical and Fourier domains.
Dutt and Rokhlin were the first to
construct an algorithm of this type, with full control of precision
\cite{nufft2}, although heuristic versions had been used earlier.
There are, by now, many variants of the NUFFT (see, for example,
\cite{nufft1,nufft2,nufft3,nufft4,nufft5,nufft6,nufft7,nufft8}).
All of these algorithms rely on
interpolation coupled with a judicious use of the FFT on an oversampled grid.
Here, we will follow the discussion in the
paper \cite{nufft6}, which describes a simple framework
for the NUFFT using Gaussian kernels for interpolation.

The type-1 NUFFT evaluates sums of the form
\begin{equation}\label{2.1}
f(\bx)=\frac{1}{N}\sum_{n=0}^{N-1}F_n e^{-i\bk_n\cdot \bx},
\end{equation}
for ``targets" $\bx$ on a regular grid  in $\bR^d$, given function values $F_n$
prescribed at arbitrary
locations $\bk_n$ in the dual space.
Here, $N$ denotes the total number of source points.

The type-2 NUFFT evaluates
sums of the form
\begin{equation}\label{2.2}
F(\bk_n)=\sum_{j_1=-M_1/2}^{M_1/2-1}\cdots \sum_{j_d=-M_d/2}^{M_d/2-1}f(\bx_j) e^{-i\bk_n \cdot \bx_j},
\end{equation}
where the ``targets" $\bk_n$ are irregularly located points in $\bR^d$,
given the function values $f(\bx_j)$ on a regular grid in the dual space.
(The type-3 NUFFT permits the sampling to be irregular in both domains, and will
not be needed in the present paper.)

We now briefly explain the basic idea underlying the NUFFT \cite{nufft2,nufft6}.
For simplicity, let us consider the one dimensional type-1 NUFFT:
\begin{equation}\label{2.3}
F(k)=\frac{1}{N}\sum_{j=0}^{N-1}f_j e^{-ikx_j}, \qquad k=-\frac{M}{2},\cdots, \frac{M}{2}-1.
\end{equation}
Note, now, that \eqref{2.3} describes the exact Fourier coefficients
of the function
\begin{equation}\label{2.4}
f(x)=\sum_{j=0}^{N-1}f_j \delta(x-x_j),
\end{equation}
viewed as a periodic function on $[0,2\pi]$. Here $\delta(x)$ denotes the Dirac function.
It is clearly not well-resolved by a uniform mesh in $x$.
By convolving with a heat kernel, however, we will construct a smooth function which
can be sampled. For this,
we let $g_\tau(x)=\sum_{l=-\infty}^{\infty}e^{-(x-2l\pi)^2/4\tau}$ denote the 1D periodic
heat kernel on $[0,2\pi]$.
If we define $f_\tau(x0$ to be convolution of $f$ with $g_\tau$:
\[ f_\tau(x)=f\ast g_\tau(x)=\int_0^{2\pi}f(y)g_\tau(x-y)dy \, ,  \]
then $f_\tau$ is a $2\pi$-periodic $C^\infty$ function and
is well-resolved by a uniform mesh in $x$ whose spacing is determined by $\tau$.
Thus, its Fourier coefficients $F_\tau(k)=\frac{1}{2\pi}\int_0^{2\pi}f_\tau(x)e^{-ikx}dx$ can be computed
with high accuracy using the standard FFT on a sufficiently fine grid.
That is,
\begin{equation}\label{2.5}
F_\tau(k)\approx \frac{1}{M_r}\sum_{m=0}^{M_r-1}f_\tau(2\pi m/M_r)e^{-ik2\pi m/M_r},
\end{equation}
where
\begin{equation}\label{2.6}
f_\tau(2\pi m/M_r)=\sum_{j=0}^{N-1}f_jg_\tau(2\pi m/M_r-x_j).
\end{equation}
Once the value $F_\tau(k)$ are known, an elementary calculation shows that
\begin{equation}
F(k)=\sqrt{\frac{\pi}{\tau}}e^{k^2\tau}F_\tau(k).
\end{equation}
This is a direct consequence of the convolution theorem and the fact that
the Fourier transform of $g_\tau$ is $G_\tau(k)=\sqrt{2\tau}e^{-k^2\tau}$.

Optimal selection of the parameters in the algorithm requires a bit of analysis,
which we omit here. We simply note \cite{nufft6} that
if $M_r=2M$ and $\tau=12/M^2$, and one uses a Gaussian to spread each
source to the nearest $24$ grid points, then the NUFFT yields about $12$ digits of accuracy.
With $\tau=6/M^2$ and
Gaussian spreading of each source to the nearest $12$ grid points, the NUFFT yields
about $6$ digits of accuracy.
The type-2 NUFFT is computed by essentially reversing the steps of type-1 NUFFT.

\section{Numerical Algorithms}

We turn now to a detailed description of our numerical algorithms for evaluating the nonlocal
 (long-range) interactions (\ref{gpe}) and the related interaction energy (\ref{gpe2}).

\subsection{High order discretization}
Since we have assumed that the function $\rho$ is smooth and rapidly decaying,
we treat it as compactly supported
with some prescribed precision $\varepsilon$ in the rectangular box
$B=[-R_1/2,R_1/2]\times\cdots\times[-R_d/2,R_d/2]$. Its Fourier
transform $\wrho$ is
\begin{equation}\label{3.1}
\wrho(\bk)=\int_Be^{- i\bk\cdot\bx}\rho(\bx)d\bx,
\end{equation}
where $\bx=(x_1,\cdots,x_d)$, $\bk=(k_1,\cdots,k_d)$.

Let us now be more specific about our
smoothness assumption. We let $\rho\in C^n(B)$, so that
$\wrho=O(\|\bk\|^{-n})$
as $\|\bk \| \rightarrow \infty$.
A straightforward calculation shows that to achieve a tolerance of  $\varepsilon$,
then evaluation of (\ref{gpe}) needs to be done only for $\|\bk\|\leq P$,
where $P=O(1/\varepsilon)^{1/n}$. We will refer to $P$ as the high-frequency cutoff.
This fixes the range of integration in $\bk$-space and bounds the oscillatory behavior of
the term $e^{- i\bk\cdot \bx}$ in the integrand of \eqref{3.1}.

Together with the fact that $\rho(\bx)$ is smooth, it follows that the
tensor product trapezoidal rule applied to
\eqref{3.1} with $N_j$ points along the $j$th axis will yield
$O(N^{-n})$ accuracy, where $N = \min_{j=1}^d N_j$. The error will decay rapidly
once each of the $N_j$ is of the order $(PR_j)$,
so that the integrand is well resolved.
If $\rho(\bx)$ is given on a uniform mesh with $N_j$ points in the $j$th dimension,
the trapezoidal rule yields
\begin{equation}\label{3.2}
\wrho(\bk)\approx \left(\prod_{j=1}^d\frac{R_j}{N_j} \right)
 \sum_{n_1=0}^{N_1-1}\cdots\sum_{n_d=0}^{N_d-1}e^{- i\bk \cdot \bx_n}
\rho(\bx_n),
\end{equation}
where $\bx_n=(-R/2+n_1(R/N_1),\cdots,-R/2+n_d(R/N_d))$.

To compute the desired solution in physical space, we need to evaluate
the inverse Fourier transform defined by (\ref{ifft}) for each of the
kernels in (\ref{int3dd}).
As discussed above, we can truncate the domain of integration in the Fourier domain
at $\|\bk\|=P=O(1/\varepsilon)^{1/n}$, with an error $\varepsilon$.
Thus, the main issue is the design of a high order rule for finite Fourier integrals
of the form:
\begin{equation}\label{3.3}
u(\bx)\approx \frac{1}{(2\pi)^d}\int_{\|\bk\|\leq P}
e^{ i\bk\cdot \bx}\wu(\bk)\wrho(\bk)d\bk \, .
\end{equation}

The principal difficulty is that the integrand above is
singular at the origin using Cartesian coordinates in $\bk$-space. It is, however,
perfectly smooth in spherical coordinates or polar coordinates, respectively.
Indeed, using the usual change of variables
in \eqref{3.3}, we obtain
\begin{equation}\label{3.5}
u(\bx)\approx \frac{1}{(2\pi)^d}
\left\{\begin{aligned}
&\int_0^P\int_0^{\pi}\int_0^{2\pi}
e^{ i\bk\cdot \bx} \|\bk\|^2\wu(\bk) \wrho(\bk)\sin\theta dk d\theta d\phi, & \text{in 3D,}\\
&\int_0^P\int_0^{2\pi}
e^{i\bk\cdot \bx} \|\bk\|\wu(\bk) \wrho(\bk)dk d\phi, &
\text{in 2D.}
\end{aligned}\right.
\end{equation}
It is easy to see that the integrand is smooth in both integrals in \eqref{3.5}
since the factor $\|\bk\|^{d-1}$ ($d=2,3$) cancels the singularity in
$\wu(\bk)$ by inspection of \eqref{int3dd}.
$\wrho(\bk)$, of course, is smooth since it is a band-limited function.

The integrals in \eqref{3.5} can be discretized with high order accuracy by using
standard (shifted and scaled) Gauss-Legendre quadrature in the radial direction (and
the longitudinal $\theta$ direction in 3D), combined with the trapezoidal rule
for the azimuthal $\phi$ variable.
Thus, we have
\begin{equation}\label{3.7}
u(\bx)\approx \frac{1}{(2\pi)^d}
\left\{\begin{aligned}
&\sum_{j_1=1}^{N_r}\sum_{j_2=1}^{N_{\theta}}\sum_{j_3=1}^{N_{\phi}}
w_je^{i\bk_j\cdot \bx}\|\bk_j\|^2 \wu(\bk_j)
\wrho(\bk_j), & \text{in 3D,}\\
&\sum_{j_1=1}^{N_r}\sum_{j_2=1}^{N_{\phi}}
w_je^{i\bk_j\cdot \bx} \|\bk_j\| \wu(\bk_j) \wrho(\bk_j), & \text{in 2D.}
\end{aligned}\right.
\end{equation}

\subsection{A simple procedure}
It is clear that $\wrho(\bk_j)$ can be evaluated from
\eqref{3.2} at the desired nonequispaced points $\bk_j$ using
the type-2 NUFFT. The summations defined in \eqref{3.7} can then be evaluated using
the type-1 NUFFT since the desired output points $\bx$ lie on a uniform grid
in physical space.

\begin{algorithm}
\caption{Simple procedure for the evaluation of \eqref{ifft}}
\label{alg1}
{\em Given the dimension $d$, the box size parameters $R_j$, $j=1,\cdots, d$
and the number of equispaced points $N_j$ in each direction,
compute $u(\bx)$ defined in \eqref{ifft} on a uniform grid in
$B=\prod_{j=1}^d[-R_j/2,R_j/2]$.}
\begin{algorithmic}[1]
\STATE Compute the coordinates on the uniform grid in $B$, that is,
$\bx_n=(-R_1/2+n_1(R_1/N_1),\cdots,-R_d/2+n_d(R_d/N_d))$, $n_j=0,\cdots, N_j$, $j=1,\cdots, d$.
\STATE Evaluate the function $\rho(\bx_n)$ at these uniform grid points.
\STATE Compute the Gauss-Legendre nodes and weights $r_j$, $w_{r_j}$, $j=1,\cdots,N_r$ for the r direction,
the trapezoidal nodes $\phi_l$, $l=1,\cdots,N_{\phi}$ for the $\phi$ direction, and
the Gauss-Legendre nodes and weights $\theta_k$, $w_{\theta_k}$, $k=1,\cdots,N_{\theta}$ for
the $\theta$ direction if $d=3$.
\STATE Use the type-2 NUFFT to evaluate $\wrho$ at these nonuniform grid points.
\STATE Use the type-1 NUFFT to evaluate $u(\bx_n)$ defined in \eqref{3.7}.
\end{algorithmic}
\end{algorithm}

The total computational cost of Algorithm \ref{alg1} is $O(N_f)+O(N_p\log N_p)$, where
$N_f$ is the total number of irregular points in the Fourier domain and $N_p$ is the total
number of equispaced points in the physical domain. As discussed in the end of Section 2,
the constant in front of $N_f$ is $24^d$ for $12$ digits accuracy with $d$ the dimension
of the problem. Since $N_f$ is often comparable with $N_p$ and the constant in the standard
FFT is quite small, the $O(N_f)$ term will dominate the computational cost in
Algorithm \ref{alg1}, making it considerably slower than the standard FFT.

\subsection{A more elaborate algorithm}
We now construct a more elaborate algorithm to reduce the interpolation cost of the
preceding scheme, by
reducing the number of irregular points $N_f$.

We note that the only singular point in the Fourier domain is the origin.
Thus we will split the integral in \eqref{3.3} into two parts using a
simple partition of unity.
\begin{equation}\label{3.9}
\begin{aligned}
u(\bx)&\approx \frac{1}{(2\pi)^d}\int_{\|\bk\|\leq P}e^{i\bk\cdot \bx}\wu(\bk)\wrho(\bk)d\bk\\
&=\int_{\|\bk\|\leq P}e^{i\bk\cdot \bx}\frac{\wu(\bk)}{(2\pi)^d}\wrho(\bk)(1-p_d(\bk))d\bk
+\int_{\|\bk\|\leq P}e^{i\bk\cdot \bx}\frac{\wu(\bk)}{(2\pi)^d}\wrho(\bk)p_d(\bk)d\bk\\
&:=I_1+I_2.
\end{aligned}
\end{equation}
We now choose the function $p_d$ so that it is
a monotone $C^\infty$ function which decays rapidly, and so that
$\frac{1-p_d(\bk)}{\|\bk\|^{d-1}}$
is a smooth function for $\bk\in\bR^d$. By the second
property of $p_d$, it is easy to see that $I_1$ can be computed using the regular FFT.
If $p_d$ decays much faster than $\hf$, then $I_2$ can be computed using the NUFFT but with
many fewer irregular points in the Fourier domain. There are many choices
for $p_d$. Indeed, any partition of unity function that is
$C^\infty$ in $\bR^d$, equals to $1$ for $\|\bk\|<R_0$ and $0$ for $\|\bk\|>R_1$ would work theoretically.
In order to minimize the oversampling factor in the evaluation of $I_1$, for the
two dimensional
problems listed in \eqref{ifft} and \eqref{int3dd}, we choose $p_2$ as follows:
\begin{equation}\label{3.12}
p_2(\bk)=\frac{1}{2}\text{erfc}\left(\frac{12(\|\bk\|-(R_0+R_1)/2)}{R_1-R_0}\right),
\end{equation}
where $erfc$ is the complementary error function defined by the formula
$\text{erfc}(x)=\frac{2}{\sqrt{\pi}}\int_x^\infty e^{-t^2}dt$.
\begin{remark}
Dimensional analysis indicates that $R_0$ and $R_1$ in the definition of
$p_2$ \eqref{3.12} should be proportional to
$O(\min(\Delta k_1,\Delta k_2))=O(1/\max(h_1,h_2))$, where $h_i$ ($i=1,2$) are the mesh size in
the $i$th direction. Our numerical experiments show that with $R_0= 0.8/max(h_1,h_2)$
and $R_1= 10/max(h_1,h_2)$, an oversampling factor of $3$ for $I_1$ and a
$60\times 60$ irregular grid
for $I_2$ yield $12$ digits accuracy.
An oversampling factor of $2$ for $I_1$ and a $40\times 40$
irregular grid for $I_2$ yield $6$ digits of accuracy.
We have not carried out a more detailed optimization.
\end{remark}

For three dimensional problems, we choose a simple Gaussian:
\begin{equation}\label{3.10}
p_3(\bk)=e^{-\frac{\|\bk\|^2}{a}},
\end{equation}
It is straightforward to verify that $\frac{1-p_3(\bk)}{\|\bk\|^2}$
has a power series expansion in $\|\bk\|^2$, so that it is a smooth function of $\bk$,
satisfies the second property.

\begin{remark}
We still need to choose the parameter $a$ in \eqref{3.10}.
Obviously, one would need fewer irregular
points for $I_2$ if $a$ were small, reducing the computational cost of $I_2$.
However, the Gaussian becomes more sharply peaked and one would need to oversample the
regular grid in $I_1$ in order to maintain high accuracy, increasing the computational cost.
Thus, $a$ should be chosen to balance the contributions of the two integral to the net cost.
Dimensional analysis indicates that $a$ should be of the order
\[ O(\min(\Delta k_1^2,\Delta k_2^2,\Delta k_3^2))=O(1/\max(h_1^2,h_2^2,h_3^2)),  \]
where the $h_i$ ($i=1,2,3$)
are the mesh size in the $i$th coordinate direction in physical space.
Numerical experiments
show that with $a=2/\max(h_1^2,h_2^2,h_3^2))$, the oversampling factor for $I_1$
can be set to $2$. An irregular (spherical)
$40\times 40\times 40$ grid in the Fourier domain achieves 12 digits of
accuracy for $I_2$, and
an irregular $24\times 24\times 24$ grid in the Fourier domain achieves 6 digits
accuracy for $I_2$.
\end{remark}

\begin{algorithm}
\caption{An improved algorithm for the evaluation of \eqref{ifft}}
\label{alg2}
{\em Given the dimension $d$, the box size parameters $R_j$, $j=1,\cdots,d$
and the number of equispaced points $N_j$ in each direction,
compute $u(\bx)$ defined in \eqref{ifft} on a uniform grid in $B=\prod_{j=1}^d[-R_j/2,R_j/2]$.}
\begin{algorithmic}[1]
\STATE Compute the coordinates of the uniform grid in $B$, that is,
$\bx_n=(-R_1/2+n_1(R_1/N_1),\cdots,-R_d/2+n_d(R_d/N_d))$, $n_j=0,\cdots, N_j$, $j=1,\cdots,d$.
\STATE Evaluate the values of the function $\rho(\bx_n)$ at these uniform grid points.
\STATE Set the oversampling factor to $2$ and compute $I_1$ in \eqref{3.9} using regular FFT.
\STATE Use NUFFT as in Algorithm \ref{alg1} to compute $I_2$ in \eqref{3.9}.
\STATE Compute $u=I_1+I_2$.
\end{algorithmic}
\end{algorithm}

\begin{remark}
Once $u(\bx)$ has been computed via Algorithm \ref{alg2},
the interaction energy \eqref{gpe2} can be
discretized via the trapezoidal rule and
evaluated by pointwise multiplication and direct summation in physical space.
The computational cost is obviously linear in the total number of discretization
points in physical space.
\end{remark}

\begin{remark}
The computational cost of the interpolation procedure within the NUFFT has been
reduced to $O(1)$ in Algorithm \ref{alg2}.
\end{remark}

\begin{remark}
Our algorithm can be easily modified to evaluate any nonlocal interaction with a
convolution structure, so long as
the Fourier transform of the kernel is known.
If the singularity at the origin of $\bk$-space cannot be
removed by switching to polar or spherical coordinates, one can
easily develop a high order generalized Gaussian quadrature rule
to discretize the singular integral in the radial direction
(see, for example, \cite{ggq2,ggq1}).
\end{remark}

\section{Numerical Examples}

We have implemented the algorithms above in Fortran.
For convenience, we have used
the publicly available software package \cite{nufft9}.
We used the gcc compiler (version $4.8.1$) with option -O3 on a $64$ bit linux
workstation with a $2.93$GHz Intel Xeon CPU and $12$ Mb of cache.
We restrict our attention to the performance of Algorithm \ref{alg2}.

For testing purposes, we have chosen right-hand sides (densities)
for which the analytical solutions are known. In the following tables,
the first column lists the total number of points in physical space, the second column
lists the value of the parameter {\it prec}.
In 2D problems $prec=0$ implies that the oversampling
factor for the regular FFT is 2 and a $40\times 40$ polar grid is used for the
NUFFT.
$prec=1$ implies that the oversampling
factor for the regular FFT is 3 and a $60\times 60$ polar grid is used for the NUFFT.
For 3D problems, $prec=0$ implies that the oversampling
factor for the regular FFT is 2 and a $24\times 24 \times 24$ spherical grid is used for
the NUFFT. $prec=1$ implies that the oversampling
factor for the regular FFT is 2 and a $40\times 40\times 40$ spherical grid is used for
the NUFFT.
The third column lists the time spent using the regular grid and the FFT.
The fourth column lists the time spent
on the irregular grid and the NUFFT. The fifth column lists the total time expended
by the algorithm. All time are measured in seconds. Finally, the last column
lists the relative $L^2$ error as compared with the analytical solution on a uniform
grid in physical space.

{\sc  Example 1: Coulomb Interactions in 2.5D.} We take $d=2$,
$U(\bx)=U_{\rm Cou}(\bx)=\frac{1}{2\pi|\bx|}$ as in (\ref{Cou2D}) and
$\rho(\bx)=e^{-|\bx|^2/a}$ with $a$ a positive constant
in (\ref{gpe}). The analytical solution to (\ref{gpe}) is given by the formula
\begin{equation}\label{4.2}
u(\bx)=\frac{\sqrt{\pi a}}{2}e^{-|\bx|^2/2a}I_0\left(\frac{|\bx|^2}{2a}\right), \qquad \bx \in \mathbb{R}^3,
\end{equation}
where $I_0$ is the modified Bessel function of order $0$ (see, for example, \cite{handbook}).
The numerical solution is computed by the formula (\ref{ifft}) with $d=2$
via Algorithm \ref{alg2}.
Table \ref{tab1} lists the numerical results for $a=1.3$.
\begin{table}[htbp]
\begin{center}
\caption{Error and timing results for Example 1.}
\begin{tabular}{|c|c|c|c|c|c|}
\hline
\emph{$N$} &\emph{Prec} &\emph{$T_{FFT}$}& \emph{$T_{NUFFT}$}& \emph{$T_{Total}$
}&\emph{$E$}\\
\hline
1024 &0 & 0.10e-2 &0.5e-2 & 0.60e-2& 1.90e-8\\
4096 &0 & 0.20e-2 &0.60e-2 & 0.80e-2& 1.57e-8\\
16438 &0 & 0.10e-1 &0.11e-1 & 0.21e-1 &1.52e-8\\
65536 &0 & 0.35e-1& 0.30e-1 &0.65e-1&1.50e-8\\
\hline
1024 &1 & 0.20e-2 &0.11e-1 & 0.13e-1& 2.40e-11\\
4096 &1 & 0.60e-2 &0.12e-1 & 0.18e-1& 3.84e-15\\
16438 &1 & 0.20e-1 &0.17e-1 & 0.37e-1 &5.75e-15\\
65536 &1 & 0.83e-1& 0.37e-1 &0.12 & 1.24e-14\\
\hline
\end{tabular}
\label{tab1}
\end{center}
\end{table}

{\sc Example 2: Dipole-Dipole Interactions with the Same Dipole Orientation in 2.5D.}
Here we take $d=2$,
$U(\bx)$ as in (\ref{dipkernel2d}) with $\alpha=0$ and $n_3=0$, and
$\rho(\bx)=e^{-|\bx|^2/a}$ with $a$ a positive constant
in (\ref{gpe}). The analytical solution to (\ref{gpe}) is given by the formula
\begin{equation}\label{4.5}
u(\bx)=\frac{3\sqrt{\pi}e^{-r}}{4\sqrt{a}}\left[I_1(r)-I_0(r)
+\frac{(\bx\cdot\bn_\perp)^2}{a}\left(2I_0(r)-\frac{1+2r}{r}I_1(r)\right)\right],\quad
\bx\in \mathbb{R}^2,
\end{equation}
where $r=\frac{|\bx|^2}{2a}$, $\bn_\perp=(n_1,n_2)^T$, and  $I_1$ is the modified Bessel
function of order $1$  (see, for example, \cite{handbook}).
The numerical solution is computed by the formula (\ref{ifft}) with $d=2$
via Algorithm \ref{alg2}.
Table \ref{tab2} lists the numerical results with $a=1.3$
and a randomly selected orientation vector $\bn_\perp=(0.52460,-0.85135)^T$.
\begin{table}[htbp]
\begin{center}
\caption{Error and timing results for Example 2.}
\begin{tabular}{|c|c|c|c|c|c|}
\hline
\emph{$N$} &\emph{Prec} &\emph{$T_{FFT}$}& \emph{$T_{NUFFT}$}& \emph{$T_{Total}$
}&\emph{$E$}\\
\hline
1024 &0 & 0.10e-2 &0.4e-2 & 0.50e-2& 1.21e-6\\
4096 &0 & 0.30e-2 &0.60e-2 & 0.90e-2& 9.33e-7\\
16438 &0 & 0.10e-1 &0.11e-1 & 0.21e-1 &8.15e-7\\
65536 &0 & 0.33e-1& 0.31e-1 &0.65e-1& 8.67e-7\\
\hline
1024 &1 & 0.30e-2 &0.10e-1 & 0.13e-1& 1.80e-8\\
4096 &1 & 0.60e-2 &0.12e-1 & 0.18e-1& 1.05e-14\\
16438 &1 & 0.20e-1 &0.17e-1 & 0.37e-1 &1.24e-14\\
65536 &1 & 0.83e-1& 0.37e-1 &0.12 & 1.74e-14\\
\hline
\end{tabular}
\label{tab2}
\end{center}
\end{table}

{\sc Example 3: Dipole-Dipole Interactions with Different Dipole Orientations in 2.5D.}
We take $d=2$,
$U(\bx)$ as  (\ref{dipkernel2dd}) with $\alpha=0$, $n_3=0$ and $m_3=0$, and
$\rho(\bx)=e^{-|\bx|^2/a}$ with $a$ a positive constant
in (\ref{gpe}). The analytical solution to (\ref{gpe}) is given by the formula
\begin{equation*}
u(\bx)=\frac{3\sqrt{\pi}e^{-r}}{4\sqrt{a}}\left[(\bn_\perp\cdot \bm_\perp)(I_1(r)-I_0(r))
+\frac{(\bx\cdot\bn_\perp)(\bx\cdot\bm_\perp)}{a}
\left(2I_0(r)-\frac{1+2r}{r}I_1(r)\right)\right],
\end{equation*}
where $r=\frac{|\bx|^2}{2a}$ and $\bm_\perp=(m_1,m_2)^T$.
The numerical solution is computed by the formula (\ref{ifft}) with $d=2$
via the Algorithm \ref{alg2}.
Table \ref{tab3} shows the numerical results with $a=1.8$, and randomly selected
orientation vectors
$\bn_\perp=(-0.44404,-0.89600)^T$ and $\bm_\perp=(0.85125,-0.52476)^T$.
\begin{table}[htbp]
\begin{center}
\caption{Error and timing results for Example 3.}
\begin{tabular}{|c|c|c|c|c|c|}
\hline
\emph{$N$} &\emph{Prec} &\emph{$T_{FFT}$}& \emph{$T_{NUFFT}$}& \emph{$T_{Total}$
}&\emph{$E$}\\
\hline
1024 &0 & 0.10e-2 &0.3e-2 & 0.40e-2& 1.43e-6\\
4096 &0 & 0.30e-2 &0.60e-2 & 0.90e-2& 1.36e-6\\
16438 &0 & 0.90e-2 &0.11e-1 & 0.20e-1 &9.99e-7\\
65536 &0 & 0.34e-1& 0.31e-1 &0.65e-1& 7.59e-7\\
\hline
1024 &1 & 0.30e-2 &0.10e-1 & 0.13e-1 & 2.40e-8\\
4096 &1 & 0.60e-2 &0.12e-1 & 0.18e-1 & 1.06e-14\\
16438 &1 & 0.21e-1 &0.17e-1 & 0.38e-1 &1.16e-14\\
65536 &1 & 0.83e-1& 0.37e-1 &0.12 & 1.88e-14\\
\hline
\end{tabular}
\label{tab3}
\end{center}
\end{table}

{\sc Example 4: Coulomb Interactions in 3D.} We take $d=3$,
$U(\bx)=U_{\rm Cou}(\bx)=\frac{1}{4\pi|\bx|}$ as (\ref{Cou3D}), and
$\rho(\bx)$ as
\begin{equation}\label{4.10}
\rho(\bx)=\left[-2\beta+4\left(\frac{x^2}{a_1^2}+\frac{y^2}{a_2^2}+\frac{z^2}{a_3^2}\right)\right]e^{-g(\bx)}, \quad g(\bx)=\frac{x^2}{a_1}+\frac{y^2}{a_2}+\frac{z^2}{a_3},
\end{equation}
where $\bx=(x,y,z)^T$, $\beta=\frac{1}{a_1}+\frac{1}{a_2}+\frac{1}{a_3}$
with $a_1$, $a_2$ and $a_3$ three positive constants,
in (\ref{gpe}). The analytical solution to (\ref{gpe}) is given by the formula
\begin{equation}
u(\bx)=-e^{-g(\bx)}=-e^{-(x^2/a_1+y^2/a_2+z^2/a_3)}, \qquad \bx=(x,y,z)^T\in \mathbb{R}^3.
\end{equation}
The numerical solution is computed by the formula (\ref{ifft}) with $d=3$
via Algorithm \ref{alg2}.
Table \ref{tab4} depicts the numerical results for $a_1=1.0$, $a_2=1.3$ and $a_3=1.5$.
\begin{table}[htbp]
\begin{center}
\caption{Error and timing results for Example 4.}
\begin{tabular}{|c|c|c|c|c|c|}
\hline
\emph{$N$} &\emph{Prec} &\emph{$T_{FFT}$}& \emph{$T_{NUFFT}$}& \emph{$T_{Total}$
}&\emph{$E$}\\
\hline
32768 &0 & 0.35e-1 &1.04 & 1.07& 2.84e-9\\
262144 &0 & 0.75 & 1.56 & 2.31& 2.51e-9\\
2097152 &0 & 5.88 & 5.06 & 10.95 &2.47e-9\\
16777216 &0 & 48.98 & 34.98 &84.07 & 2.46e-9\\
\hline
32768 &1 & 0.36e-1 &4.65 & 4.69& 9.06e-10\\
262144 &1 & 0.75 & 5.70 & 6.45& 8.72e-14\\
2097152 &1 & 5.93 & 9.42 & 15.36 & 8.53e-14\\
16777216 &1 & 49.92 & 49.98 & 100.00 & 8.51e-14\\
\hline
\end{tabular}
\label{tab4}
\end{center}
\end{table}

{\sc Example 5: Dipole-Dipole Interactions with the Same Dipole Orientation in 3D.}
We take $d=3$, with
$U(\bx)$ given by (\ref{dipkernel2d}) and
$\rho(\bx)$ given by  (\ref{4.10}) in (\ref{gpe}).
The analytical solution to (\ref{gpe}) is given by the formula
\begin{equation*}
u(\bx)=-\rho(\bx)+6\left[2g_\bn(\bx)^2
-\left(\frac{n_1^2}{a_1}+\frac{n_2^2}{a_2}+\frac{n_3^2}{a_3}\right)\right]e^{-g(\bx)},
\quad g_\bn(\bx)=\frac{x n_1}{a_1}+\frac{y n_2}{a_2}+\frac{z n_3}{a_3}.
\end{equation*}
The numerical solution is computed by the formula (\ref{ifft}) with $d=3$
via the Algorithm \ref{alg2}.
Table \ref{tab5} lists the numerical results for
 $a_1=1.3$, $a_2=1.5$, $a_3=1.8$ and randomly selected orientation vector
 $\bn=(-0.36589,-0.69481,0.61916)^T$.
\begin{table}[htbp]
\begin{center}
\caption{Error and timing results for Example 5.}
\begin{tabular}{|c|c|c|c|c|c|}
\hline
\emph{$N$} &\emph{Prec} &\emph{$T_{FFT}$}& \emph{$T_{NUFFT}$}& \emph{$T_{Total}$
}&\emph{$E$}\\
\hline
32768 &0 & 0.36e-1 &1.04 & 1.08& 3.06e-7\\
262144 &0 & 0.75 & 1.55 & 2.30& 1.15e-7\\
2097152 &0 & 5.84 & 5.07 & 10.93 &3.81e-8\\
16777216 &0 & 50.53 & 35.77 &86.39 & 1.13e-7\\
\hline
32768 &1 & 0.37e-1 &4.65 & 4.69& 2.65e-7\\
262144 &1 & 0.75 & 5.70 & 6.44& 1.13e-13\\
2097152 &1 & 5.90 & 9.39 & 15.30 & 1.12e-13\\
16777216 &1 & 50.53 & 49.38 & 100.01 & 8.67e-14\\
\hline
\end{tabular}
\label{tab5}
\end{center}
\end{table}

{\sc Example 6: Dipole-dipole Interactions with Different Dipole Orientations in 3D.}
We take $d=3$, with
$U(\bx)$ given by (\ref{dipkernel2dd}) with $\alpha=0$,  and
$\rho(\bx)$ given by (\ref{4.10}) in (\ref{gpe}).
The analytical solution to (\ref{gpe}) is given by
the formula
\begin{equation}\label{4.15}
u(\bx)=-\rho(\bx)+6
\left[2g_\bn(\bx)g_\bm(\bx)
-\left(\frac{n_1 m_1}{a_1}+\frac{n_2 m_2}{a_2}+\frac{n_3 m_3}{a_3}\right)\right]  e^{-g(\bx)},
\quad \bx\in\mathbb{R}^3.
\end{equation}
The numerical solution is computed by the formula (\ref{ifft}) with $d=3$
via Algorithm \ref{alg2}.
Table \ref{tab6} shows the numerical results for
$a_1=1.2$, $a_2=1.45$, $a_3=1.73$, $\bn=(0.82778,0.41505,-0.37751)^T$,
$\bm=(0.31180,0.93780,-0.15214)^T$.
\begin{table}[htbp]
\begin{center}
\caption{Error and timing results for Example 6.}
\begin{tabular}{|c|c|c|c|c|c|}
\hline
\emph{$N$} &\emph{Prec} &\emph{$T_{FFT}$}& \emph{$T_{NUFFT}$}& \emph{$T_{Total}$
}&\emph{$E$}\\
\hline
32768 &0 & 0.36e-1 &1.04 & 1.08& 3.18e-7\\
262144 &0 & 0.75 & 1.55 & 2.30& 1.59e-7\\
2097152 &0 & 5.92 & 5.11 & 11.04 &8.28e-8\\
16777216 &0 & 50.20 & 35.80  &86.09  & 1.53e-7\\
\hline
32768 &1 & 0.36e-1 &4.65 & 4.69& 3.36e-7\\
262144 &1 & 0.75 & 5.70 & 6.46& 1.00e-13\\
2097152 &1 & 5.88 & 9.36 & 15.26 & 3.90e-13\\
16777216 &1 & 50.34 & 42.36 & 92.79 & 1.04e-13\\
\hline
\end{tabular}
\label{tab6}
\end{center}
\end{table}

{\sc Example 7: The Interaction Energy of Dipole-dipole Interactions with the Same Dipole
Orientation in 3D.}
We take $d=3$, with $U(\bx)$ given by (\ref{dipkernel2d}) $\bn=(0, 0, 1)^T$ and
$\rho(\bx)$ given by
\begin{equation}
\rho(\bx)=\pi^{-3/2}\gamma_x \sqrt{\gamma_z} e^{-(\gamma_x(x^2+y^2)+\gamma_zz^2)},
\qquad \bx=(x,y,z)^T\in \mathbb{R}^3,
\end{equation}
with $\gamma_x$, $\gamma_y$ and $\gamma_z$ three positive constants,  in (\ref{gpe}).
The dipole-dipole interaction energy $E(\rho)$ in (\ref{gpe2}) can be
evaluated analytically as \cite{bao1,Parker,Tick,Tikhonenkov}
\begin{equation}
E(\rho)=-\frac{\lambda \gamma_x \sqrt{\gamma_z}}{4\pi \sqrt{2\pi}}
\left\{\begin{aligned}
&\frac{1+2\kappa^2}{1-\kappa^2}-\frac{3\kappa^2\arctan\sqrt{\kappa^2-1}}
{(1-\kappa^2)\sqrt{\kappa^2-1}}, & \kappa>1,\\
&0, & \kappa=1,\\
&\frac{1+2\kappa^2}{1-\kappa^2}-\frac{3\kappa^2}{2(1-\kappa^2)\sqrt{1-\kappa^2}}
\ln\left(\frac{1+\sqrt{1-\kappa^2}}{1-\sqrt{1-\kappa^2}}\right), & \kappa<1, \\
\end{aligned}\right.
\end{equation}
where $\kappa=\sqrt{\gamma_z/\gamma_x}$. Three cases with $\lambda =8\pi/3$ are considered here.

\smallskip

Case I. $\gamma_x=0.25$ and $\gamma_z=1$, the exact energy is $E(\rho)\approx 0.03867086140999021$;

Case II. $\gamma_x=1$ and  $\gamma_z=1$, the exact energy is
 $E(\rho)= 0$;

Case III. $\gamma_x=2$ and $\gamma_z=1$, the exact energy is $E(\rho)\approx -0.1386449740987819$.

\smallskip

The energy is computed numerically by the formula (\ref{ifft}) with $d=3$
via Algorithm \ref{alg2} (see Remark 3.3 as well).
Table \ref{tab7} shows the numerical results for the above three cases.
From this table, we observe that our algorithm achieves full double precision,
while at most seven-digit accuracy is achieved in \cite{bao1}.

\begin{table}
\begin{center}
\caption{Numerical results for computing the dipole-dipole interaction energy in 3D.
Here $N$ is the total number of points in the computational domain, $N_I$ is the total
number of irregular points in the Fourier domain, $E_c$ is the computed value of the
dipole-dipole interaction energy, $E:=|E_c-E(\rho)|$ is the absolute error, and
$T$ is the total CPU time in seconds.}
\center\begin{tabular}{|c|c|c|c|c|c|}
\hline
 & $N$ & $N_I$ & $E_c$& $E$ & T \\
\hline
& 32768 & 1000   &0.03867932878216508 & 8.5e-6  & 0.14\\
& 32768 & 3375   &0.03867085889326449 & 2.5e-9  & 0.25\\
Case I & 262144 & 8000 &0.03867086140931093 & 6.8e-13  & 1.12\\
& 262144 & 27000 &0.03867086140998955 & 6.7e-16 & 1.89\\
\hline
& 32768 & 1000    & -4.096001975539615e-7 & 4.1e-7  & 0.14\\
& 32768 & 3375   & 1.983944760166238e-9 & 2.0e-9  & 0.23\\
Case II & 262144 & 8000 & 5.117228774684975e-14 & 5.1e-14 & 1.19\\
& 262144 & 27000 & -7.785768552598901e-16  & 7.8e-16 & 1.87\\
\hline
& 32768 & 1000   &$-0.1386463326990541$ & 1.4e-6  & 0.14\\
& 32768 & 3375   &$-0.1386449725315638$ & 1.6e-9  & 0.20\\
Case III & 262144 & 8000 &$-0.1386449740987009$ & 8.1e-14  & 1.12\\
& 262144 & 27000 &$-0.1386449740987584$ & 2.3e-14 & 1.84\\
\hline
\end{tabular}
\label{tab7}
\end{center}
\end{table}

\begin{remark}
From these tables, it is clear that the timing scales roughly linearly
with $N$. The timing difference between $prec=0$
and $prec=1$ are not very significant. Thus, we recommend setting $prec=1$
in general.
\end{remark}

\section{Conclusions}
 An efficient and high-order algorithm has been constructed for the evaluation of
long-range Coulomb and dipole-dipole interactions of the type which arise
in quantum physics and chemistry, as well as materials
simulation and design. The algorithm evaluates these interactions
in the Fourier domain, with a coordinate transformation that removes the
singularity at the origin. The Fourier integral is then
discretized via high-order accurate quadrature,
and the resulting discrete summation is carried out using the nonuniform FFT (NUFFT).
The algorithm is straightforward
to implement and requires $O(N\log N)$ work, where $N$ is the total number
of points in the physical space discretization.
Thus, the net cost is of the same order as that of using the uniform FFT for
problems with periodic boundary conditions.
Our algorithm is easily extended to the computation of any nonlocal interaction
that has a suitable convolution structure.
When the singularity in the Fourier transform cannot be accounted for by a simple change
of variables, generalized Gaussian
quadrature can be used to create a  high-order discretization \cite{ggq2,ggq1}, to
which the NUFFT can be applied, achieving  nearly optimal computational complexity.
An example where this extension of our algorithm is helpful is the solution
of the free space Poisson equation in 2D, where
the Fourier transform of the kernel is $\frac{1}{\|\bk\|^2}$
and the integrand is still singular after switching to polar coordinates.
We plan to incorporate the method described here
into efficient and accurate solvers for computing
the ground state and dynamics of dipolar BECs, the nonlinear Sch\"{o}dinger equation
with a Coulomb potential, and the Kohn-Sham equations for electronic structure.
Similar ideas have been used for computing Stokes interactions with compactly
supported data with a mixture of free space and periodic boundary conditions
imposed on a unit cell \cite{lindbo,lindbo2}.

\vspace{.4in}

\section*{Acknowledgments}
The authors would like to thank Mark Tygert and Anna-Karin Tornberg
for many useful discussions.
Part of this work was done when the first and third authors were visiting
Beijing Computational Science Research Center in the summer of 2013.

\end{document}